\documentclass[12pt]{article}
\usepackage{amsmath,amsthm,amssymb,latexsym,amscd}
\usepackage{setspace}
\usepackage[dvips]{graphicx} 
                                                                         
\usepackage[dvips]{color} 
\newtheorem{theorem}{Theorem}[section]
\newtheorem{definition}{Definition}
\newtheorem{lemma}[theorem]{Lemma}
\newtheorem{corollary}[theorem]{Corollary}
\newtheorem{proposition}[theorem]{Proposition}

\newcommand{\R}{\mathbb{R}}
\newcommand{\Z}{\mathbb{Z}}
\newcommand{\vv}{{\bf{v}}}

\newcommand{\mc}{\mathcal{F}}
\newcommand{\mco}{\mathcal{F}_1}
\newcommand{\mct}{\mathcal{F}_2}
\newcommand{\mch}{\mathcal{F}_3}

\newcommand{\mlink}{(\mco,\mct,\mch)}
\newcommand{\mpfot}{M(\mco,\mct)}

\newcommand{\mpfoth}{M(\mco,\mct,\mch)}
\newcommand{\fpr}{M(\mco,\mct)\times_{M(\mco)} M(\mco,\mch)}
\newcommand{\vev}{(v_1 , v_2 , v_3 , v_4 )}

\newcommand{\thio}{\theta_{i,1}}
\newcommand{\thnoo}{\theta_{{n_1},1}}
\newcommand{\thot}{\theta_{1,2}}
\newcommand{\thit}{\theta_{i,2}}
\newcommand{\thntt}{\theta_{{n_2},2}}
\newcommand{\thto}{\theta_{2,1}}

\newcommand{\thoh}{\theta_{1,3}}

\newcommand{\thih}{\theta_{i,3}}

\newcommand{\thnhh}{\theta_{{n_3},3}}

\newcommand{\prest}{\left. P \right|_{\backslash \mathcal{F}}}

\newcommand{\prestoh}{\left. P_1 \right|_{\backslash \mathcal{F}_3}}
\newcommand{\presttt}{\left. P_2 \right|_{\backslash \mathcal{F}_2}}
\begin{document}

\title{On the Moduli Space of Multipolygonal Linkages in the Plane}

\author{Michael Holcomb\\
Mathematics Department\\
Louisiana State University\\
Baton Rouge, Louisiana\\
holcomb@math.lsu.edu}

\date{6/24/03}

\maketitle

\begin{abstract}

The geometric, topological, and symplectic properties of moduli spaces (spaces of configurations modulo rotations and translations) of polygonal linkages have been studied by Kapovich, Millson, and Kamiyama, et. al.  One can form a polygonal linkage by taking two free linkages and identifying initial and terminal vertices.  This can be generalized so that one takes three free linkages and identifies initial and terminal vertices.  Then one obtains a linkage which contains multiple polygons, any two of which have shared edges.  The geometric and topological properties of moduli spaces of these multipolygonal linkages are studied.  These spaces turn out to be compact algebraic varieties.  Some conditions under which these spaces are smooth manifolds, cross products or disjoint unions of moduli spaces of polygonal linkages, or connected, are determined.  In addition, dimensions in the smooth manifold cases and some Euler characteristics are computed.
\end{abstract}

\section{Introduction}
\label{S:introduction}

In robotics one is interested in the totality of possible configurations of a robot arm or, more generally, a mechanical linkage.  One starts with a collection of bars, or edges, and fastens them together at hinges, or vertices, that allow either two-dimensional or three-dimensional mobility.  In a robot arm, which is a free linkage, we can fix an initial edge, make each hinge composed of exactly two consecutive bars, and let the terminal vertex move freely.  The space of possible configurations is just a product of circles.

In more general linkages, we can fix more vertices, make any hinge connect two or more edges, and some vertices can be identified.  For example, if we take two free linkages and identify their initial and terminal vertices, then we obtain a polygonal linkage.

There is a substantial amount of literature on the topological properties of moduli spaces of special polygonal linkages, particularly equilateral polygonal linkages.  A polygonal linkage in $\R^n$ is simply a polygon whose vertices are hinges and whose edges are bars.   We are allowed to deform and collapse the polygon without changing the lengths of any bars or any incidences of vertices and edges.  This process can change some (or possibly all) of the angles at hinges.  An equilateral polygonal linkage is one where each bar has the same length.

The configuration space of a polygonal linkage is the space of realizations of the linkage.  The generic dimension of this space depends on the number of edges in the linkage and the ambient space.  For example, the dimension of the moduli space of a generic hexagon in $\R^3$ is 18.  But it is convenient to reduce this dimension by fixing one vertex at the origin and one adjacent vertex along one of the coordinate axes.  Mathematically speaking, we take space of orbits of configurations for the group of Euclidean motions (translations and rotations).  We call this orbit space the moduli space for the linkage.  

Even though these moduli spaces cannot be visualized in general, geometric properties such as dimension, smoothness, Betti numbers, etc., can be investigated.  Let $M({d_1},\ldots,{d_n})$ be the moduli space of polygonal linkages in $\R^2$, with side lengths $d_i$, normalized so that the perimeter ${d_1}+\ldots+{d_n}$ equals 1.  In order for $M$ to be nonempty, the $d_i$ must belong to a polytope $P$ in $\R^2$.  Kapovich and Millson showed that the topological type of $M$ depends on the vector ${\bf{d}}=({d_1},\ldots,{d_n})$ in the following way.  There is a subdivision of the polytope $P$ into chambers,  The chambers are divided by walls such that the topological type of $M({\bf{d}})$ is constant for ${\bf{d}}$ in each chamber but undergoes a change as ${\bf{d}}$ crosses a wall.  Kamiyama studied Euler characteristics in \cite{Kte} and \cite{Kec}, homology groups and Poincar\'{e} Polynomials in \cite{Kho} and \cite{Kmo}, Chern numbers in \cite{Kcn}, symplectic volume in \cite{KTs}, and conditions in which the moduli space is a smooth manifold in \cite{Kte}.  

In this thesis we explore combinations of two polygonal linkages which share at least one edge.  We build polygonal and multipolygonal linkages up from free linkages.  The notation for these linkages is set up in Section \ref{S:polyfree} and \ref{S:multpoly}.  The moduli spaces are compact, Hausdorff, real triangulable algebraic varieties.  In Section \ref{S:multiquad} we investigate multiquadrilateral linkages and classify the moduli spaces which have dimension less than two.  

In Section \ref{S:gen} we find topological results for certain classes of multipolygonal linkages.  First we find two classes of multipolygonal linkages that are smooth manifolds.  Next, we observe that moduli spaces of multipolygonal linkages can be thought of as a fibered product of moduli spaces of polygonal linkages.  Then we join two ``long enough'' edges to a polygonal linkage and determine when the resulting moduli spaces is either equal to or is close to a disjoint union of copies of the original polygonal moduli space.  Finally, we find a class of multipolygonal linkages whose moduli space is connected.

In Section \ref{S:spec} we map the moduli space to $S^1$ via $p$, by sending a configuration to a pre-specified angle in that configuration.  We make use of some technical theorems in differential geometry and algebraic topology to determine when all fibers of the map $p$ are homeomorphic, and what can cause the fibers to change topological type and at what angle values.  We conclude the section with a specific example involving an Euler characteristic computation.

\section{Polygonal Linkages and Free Linkages}
\label{S:polyfree}

We define a free linkage in $\mathbb{R}^p$ in terms of the lengths of each edge.  Let $${\bf{d}}=({d_1},\ldots,{d_n})\in \mathbb{R}_{>0}^p$$ be a vector of positive lengths.  The set
$$\mathcal{F}(p;{\bf{d}})=\{\left. (x_1,\ldots,x_{n+1})\in{\R^p}^n \right| |{x_i}-{x_{i+1}}|=d_i \quad \text{for} \quad i=1,\ldots,n\}$$ 
denotes the space of embeddings of linkages with length vector ${\bf{d}}$.  Since we will study linkages in the plane, we will write $\mathcal{F}(\bf{d})$ for $\mathcal{F}(2;\bf{d})$.  When the lengths are understood, we will simply use $\mc$.

Note that we will use $\mc$ as a particular embedding, but we will also use $\mathcal{F}$ as a generic object.  Sometimes we will say ``moving continuously from one point to another point in $\mc$'', and other times we will say ``continuously deforming $\mc$''. The difference in usage should be clear from the context.  

We define $l_{\mathcal{F}}$ to be the length between the initial vertex $x_1$ and the terminal vertex $x_{n+1}$ in $\mc$.

When length vectors have many repeated edges, we will let $d^{<m>}$ denote repeated lengths.  For example, we will write $(2,2,1,1,1,1)$ as $(2^{<2>};1^{<4>})$.  And if $\bf{d}$ has $n$ edges of length $a$, then we simply write ${\bf{d}} = (a^{<n>})$.

We can now form more linkages composed of free linkages.  We will use the following notation for two or more ordered free linkages that have the same initial and terminal points.
$$(\mathcal{F}_1,\ldots, \mathcal{F}_r), \quad{\mathcal{F}_j} \in{F(p;{\bf{d}_j})}$$

Notice that $(\mathcal{F}_1,\mathcal{F}_2)$ denotes a polygonal linkage.  For example, if we specify a priori a free linkage $\mco(3,4)$ and a free linkage $\mct(2,2,1,3)$, and we further force the initial and terminal points of $\mco$ and $\mct$ to agree, then we obtain a hexagonal linkage with lengths in the following order:  3,4,3,1,2,2.

\begin{figure}[ht]
\scalebox{0.60}{ \includegraphics{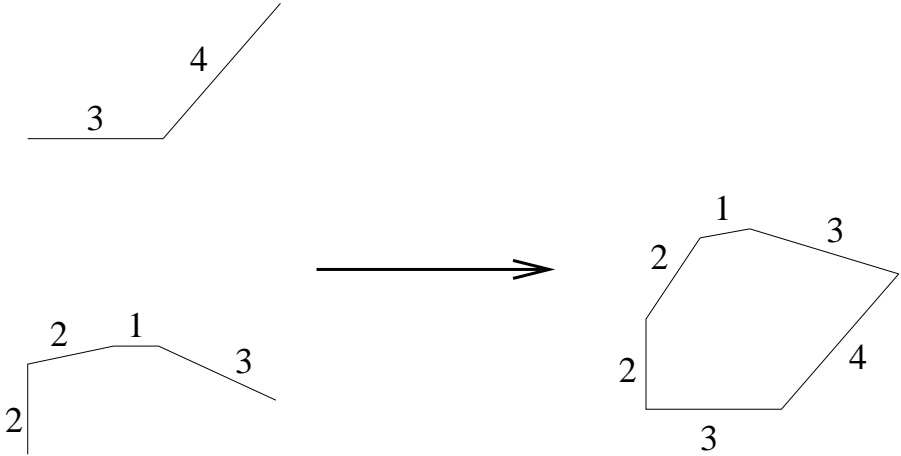}}
\caption{Forming a Hexagonal Linkage from Two Free Linkages}
\label{Fi:hexintwo}
\end{figure}

In this paper we will study multipolygonal linkages composed of three free linkages embedded in $\mathbb{R}^2$.  They will be denoted $(\mathcal{F}_1,\mathcal{F}_2,\mathcal{F}_3)$.  We call them multipolygonal because any pair of free linkages forms a polygonal linkage. 

We use $C(\mathcal{F}_1,\mathcal{F}_2,\mathcal{F}_3)$ to denote the space of all embeddings of $\mco$, $\mct$, and $\mch$ with the initial points identified at the origin and the terminal points identified.  This space has the following nice properties.

\begin{proposition}
\label{P:niceprops}
The space $C(\mathcal{F}_1,\mathcal{F}_2,\mathcal{F}_3)$ is a compact, Hausdorff, real algebraic variety.
\end{proposition}

More generally, the space $C(\mathcal{F}_1,\ldots,\mathcal{F}_r)$ is defined similarly and has the same nice properties.

For example, if $\mc$ has $n$ edges, then $C(\mc)$ equals a product of $n$ circles.

In addition, we can mod out by the group $SO(2)$ of orientation-preserving isometries of ${\mathbb{R}}^2$.  Because we already specified that initial and terminal points of each free linkage must agree, we can rotate the first edge of only one free linkage (we will do this to $\mathcal{F}_1$) so that it lies on the positive $x$-axis.  

Let $$M(\mathcal{F}_1,\mathcal{F}_2,\mathcal{F}_3)=C(\mathcal{F}_1,\mathcal{F}_2,\mathcal{F}_3)/SO(2).$$
These denote moduli spaces, which are sets of equivalence classes of embeddings of $(\mco,\mct,\mch)$ in $\R^2$ modulo the action of the special Euclidean group $SE(2,\R)=\R^2 \rtimes SO(2)$.

We have the following propositions about $M(\mathcal{F}_1,\mathcal{F}_2,\mathcal{F}_3)$.

\begin{proposition}
\label{P:niceprops2}
The space $M(\mathcal{F}_1,\mathcal{F}_2,\mathcal{F}_3)$ is a compact, Hausdorff, real algebraic variety.
\end{proposition}

\begin{proposition}
\label{P:triangulable}
The space $M(\mathcal{F}_1,\mathcal{F}_2,\mathcal{F}_3)$ is triangulable.
\end{proposition}

See Th\'{e}or\`{e}me 9.2.1 in \cite{Tri} for the proof of the second proposition.  More generally, the space $M(\mathcal{F}_1,\ldots,\mathcal{F}_r)$ is defined similarly and has the same propositions true about them.

We will study moduli spaces of the form $M(\mathcal{F}_1,\mathcal{F}_2,\mathcal{F}_3)$.  Points in the moduli space will be represented by $P$ with possible subscripts.  

If we rearrange the free linkages by a nontrivial permutation $\sigma$, then we obtain a space $M(\mathcal{F}_{\sigma(1)},\ldots,\mathcal{F}_{\sigma(n)})$ which is different from, but homeomorphic to, the space $M(\mathcal{F}_1,\ldots,\mathcal{F}_n)$.

\section{Multipolygonal Linkages}
\label{S:multpoly}

We use the following notation with a multipolygonal linkage $\mlink$ in the plane.  The length vector for $\mathcal{F}_i$ will equal $(d_{1,i},\ldots,d_{n,i})$.  We will denote the $j$th vertex of the $i$th free linkage by $x_{j,i}$.  The edge $e_{j,i}$ will link $x_{j,i}$ to $x_{j+1,i}$.  When we reference the angle between $e_{j,i}$ and the $x$-axis, we will think of $e_{j,i}$ as the vector $d_{j,i}e^{\sqrt{-1}\theta_{j,i}}$ directed from $x_{j,i}$ to $x_{j+1,i}$.

We identify $x_{1,1}$, $x_{1,2}$, and $x_{1,3}$, the initial points of the free linkages, at the origin, and we identify $x_{{n_1},1}$, $x_{{n_2},2}$, and $x_{{n_3},3}$.

\begin{figure}[ht]
\scalebox{0.60}{ \includegraphics{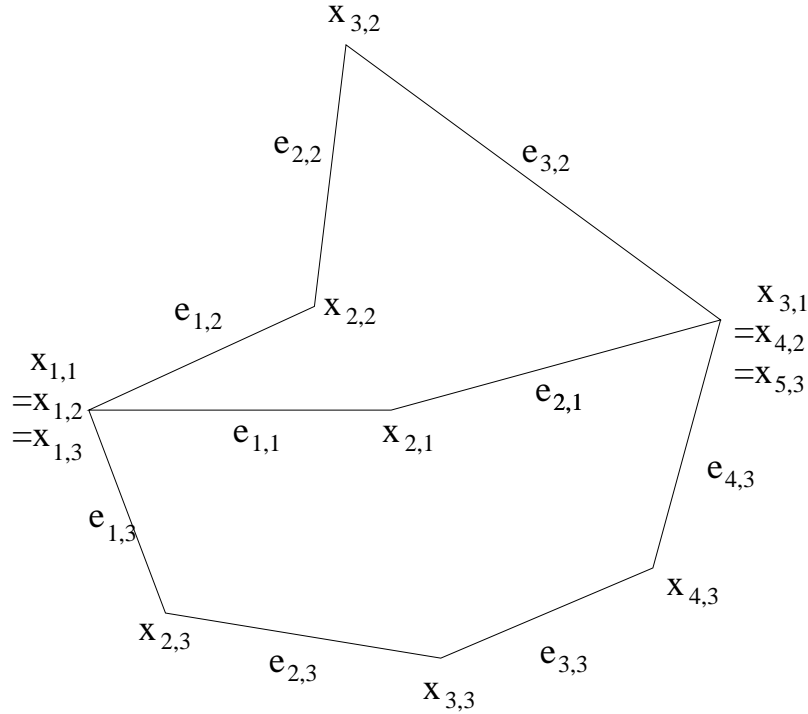}}
\caption{A Generic Multipolgonal Linkage Formed by Three Free Linkages}
\label{Fi:generic}
\end{figure}

\section{Multiquadrilateral Linkages}
\label{S:multiquad}

If ${n_1}={n_2}={n_3}=2$, then we have multiple joined quadrilateral linkages.  To investigate the moduli spaces of these linkages, it will be convenient to investigate the sublinkages which contain $\mco\in \mathcal{F}(d_{1,1},d_{2,1})$.  So we will denote the linkage which contains only $\mct \in \mathcal{F}(d_{1,2},d_{2,2})$ and $\mco$ by $(\mco,\mct)$.  Likewise, we will denote the linkage which contains only $\mch\in \mathcal{F}(d_{1,3},d_{2,3})$ and $\mco$ by $(\mco,\mch)$.  

For simple notation, the lengths $a$, $b$, $c$, $d$, $e$, and $f$ will refer to $d_{1,1}$, $d_{2,1}$, $d_{1,2}$, $d_{2,2}$, $d_{1,3}$, and $d_{2,3}$ respectively.  We can assume that 
\begin{equation}
\label{E:minpair}
min\{a+b,c+d,e+f\} = a+b
\end{equation}
so that $\mco$ has sum of lengths of its edges minimal.  Then $l$ refers to the common length from $x_{1,1}$ to $x_{3,1}$, from $x_{1,2}$ to $x_{3,2}$, and from $x_{1,3}$ to $x_{3,3}$.  Finally, the line containing $x_{1,1}$ and $x_{3,1}$ will be referred to as $L$.

\begin{figure}[ht]
\scalebox{0.60}{ \includegraphics{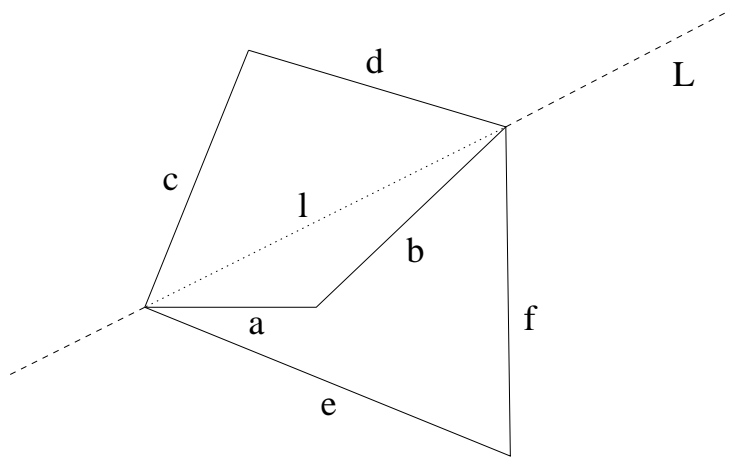}}
\caption{A generic multiquadrilateral linkage}
\label{Fi:quadgeneric}
\end{figure}

\begin{definition}
\label{D:generic}
A multiquadrilateral linkage $(\mco,\mct,\mch)$ is called generic if $a=b$, $c=d$, and $e=f$ are not simultaneously true.
\end{definition}

We will only investigate the generic case.  The moduli spaces of these generic multiquadrilateral linkages will have dimension less than two. 

For each $i$, let $\theta_i$ be the clockwise angle from $e_{1,i}$ to $e_{2,i}$.  We focus mostly on the angle $\theta_1$.  The map $p:\mpfoth\rightarrow S^1$ sends a particular configuration $(\mco,\mct,\mch)$ to the angle $\theta_1$.  The lengths of the edges in $(\mco,\mct)$ can restrict the angle $\theta_1$.  Likewise, the lengths of the edges in $(\mco,\mch)$ can restrict $\theta_1$.  So we have restriction maps $\left. p\right|_{(\mco,\mct)}$ whose image is $\Theta_1$ and $\left. p\right|_{(\mco,\mch)}$ whose image is $\Theta_2$.  Each $\Theta_i$ is a subset of $[0,2\pi)$.

\begin{lemma}
\label{L:restrict}
The image of the map $p$ is $\Theta_1 \cap \Theta_2$.
\end{lemma}

Because of condition (\ref{E:minpair}), the image of $p$ is either empty, the angle $\pi$, a closed interval containing $\pi$, or all of $[0,2 \pi)$.

\begin{definition}
\label{D:collapse}
The free linkage $\mathcal{F}_i$ is collapsible if there exists at least one configuration of $(\mco,\mct,\mch)$ where $\mathcal{F}_i$ lies completely in a straight line such that $\theta_i = \pi$.
\end{definition}

\begin{definition}
\label{D:fold}
The free linkage $\mathcal{F}_i$ is foldable if there exists at least one configuration of $(\mco,\mct,\mch)$ where $\mathcal{F}_i$ lies completely in a straight line such that $\theta_i = 0$. 
\end{definition}

Note that $\mco$ is already collapsible because of condition (\ref{E:minpair}).  Here are the conditions on the lengths that determine if $\mathcal{F}_i$ is collapsible or foldable.

\begin{align}
&\mct \text{\qquad collapsible:\qquad} && a+b=c+d.\notag\\
&\mct \text{\qquad foldable:\qquad}    && |c-d|\geq|a-b| \text{\quad and \quad} |c-d|\geq|e-f|.\notag\\
&\mch \text{\qquad collapsible:\qquad} && a+b=e+f.\notag\\
&\mch \text{\qquad foldable:\qquad}    && |e-f|\geq|a-b| \text{\quad and \quad} |e-f|\geq|c-d|.\notag
\end{align}
Here is an example.  The length constraints 
$$a+b<c+d$$
$$a-b<c-d$$
on $(\mco,\mct)$ impose the angle constraint $\theta_{min}\leq{\theta_1}\leq{\theta_{max}}$, where 
$$\theta_{min} = cos^{-1}\frac{a^2+b^2-(d-c)^2}{2ab}$$ and
$$\theta_{max} = 2\pi - \theta_{min}.$$

Here are the possible subsets $im(p)$ described above.

\begin{enumerate}
\item If $im(p)$ is empty, then it is because $a+b<|c-d|$ or $a+b<|e-f|$.  In this case, no angles are possible for $\theta_1$, and the moduli space is empty.
\item If $im(p) = \{\pi\}$, then it is because $a+b = |c-d|$ and $a+b\geq{|e-f|}$, or because $a+b=|e-f|$ and $a+b\geq{|c-d|}$. In this case, the only possible angle for $\theta_1$ is $\pi$ and the moduli space is either one or two points, with one point occuring if the inequality is actually an equality.
\item If $im(p)$ is an interval containing $0$, then $im(p)$ is all of $[0,2\pi)$.  Here the moduli space is four circles with possible intersections. (CASE$(\mathbf{A})$)
\item If $im(p)$ is an interval not containing $0$, then $im(p)$ is of the form $[\theta_{min},\theta_{max}]\subset(0,2\pi)$, and the moduli space is one circle or a figure-8.(CASE $(\mathbf{B})$)
\end{enumerate}

CASE $\bf{(A)}$:  Both $0$ and $\pi$ are possible angles for $\theta_1$.

In terms of collapsibility and foldability, this case means that $\mco$ is collapsible and foldable.  We define a vector ${\bf{v}}=\vev$ of zeros and ones as follows.  We let 

\[
v_1 =
\begin{cases}
1 &\text{if $\mct$ is collapsible;}\\
0 &\text{if $\mct$ is not collapsible,}
\end{cases} 
\]

\[
v_2 =
\begin{cases}
1 &\text{if $\mct$ is foldable;}\\
0 &\text{if $\mct$ is not foldable,}
\end{cases} 
\]

\[
v_3 =
\begin{cases}
1 &\text{if $\mch$ is collapsible;}\\
0 &\text{if $\mch$ is not collapsible,}
\end{cases} 
\]
and
\[
v_4 =
\begin{cases}
1 &\text{if $\mch$ is foldable;}\\
0 &\text{if $\mch$ is not foldable.}
\end{cases} 
\]
We also define the pairs of ${\bf{v}}$ as $({v_1},{v_2})$ and $({v_3},{v_4})$.

\begin{theorem}
\label{T:fourcircles}
 (Case $\bf{(A)}$)  Suppose that $(\mco,\mct,\mch)$ is generic and that $0$ and $\pi$ are both possible angles for $\theta_1$.  Then the moduli space of the multiquadrilateral linkage $\mpfoth$ is obtained from four circles, with possible identifications.  These can be determined as follows. 
\begin{align}
&\text{Number of components\qquad} &&=2^{\text{\# of (0,0) pairs in $\vv$}}\notag\\
&\text{Points in $p^{-1}(0)$\qquad} &&=2^{2-\text{\# of 1's in even positions of $\vv$}}\notag\\
&\text{Points in $p^{-1}(\pi)$\qquad}&&=2^{2-\text{\# of 1's in odd positions of $\vv$}}\notag
\end{align}
\end{theorem}

\begin{proof} We prove the special case of $v=(0,0,1,0)$, and the other cases follow a similar proof. 

We let $P(\theta_2 , \theta_3 )$  be any point in $p^{-1}(\theta_1)$, where initially $\theta_1 \neq{0,\pi}$ and where $\theta_2$ and $\theta_3$ depend on $\theta_1$.  Suppose $\theta_2 \in (0,\pi)$ and $\theta_3 \in (\pi,2 \pi)$ initially.  We denote a point in $p^{-1}(0)$ by $P(\theta_2 (0),\theta_3 (0))$, and we denote a point in $p^{-1}(\pi)$ by $P(\theta_2 (\pi),\theta_3 (\pi))$.

First, we find a circle of configurations starting at a point $P(\theta_2 , \theta_3 )$.  We can perturb $\theta_1$ continuously on the interval $[0,2\pi)$.  As we do so, the angle $\theta_2$ remains in $(0,\pi)$ since $\mct$ is neither collapsible nor foldable.  But the angle $\theta_3$ remains in $[\pi,2 \pi)$ since $\mch$ is collapsible but not foldable.  

The path-component containing each configuration in $p^{-1}(\theta_1)$ contains at least that circle of configurations.  Whether the component is bigger depends on whether any other circles of configurations meet it when $\theta_1$ equals $0$ or $\pi$.  This in turn depends on how many points lie in $p^{-1}(0)$ and $p^{-1}(\pi)$. 

There are four points in $p^{-1}(\theta_1)$.

\begin{equation}
P_1 = P(\theta_2 ,\theta_3 )\notag
\end{equation}
\begin{equation}
P_2 = P(\theta_2 ,2\pi-\theta_3 )\notag
\end{equation}
\begin{equation}
P_3 = P(2\pi-\theta_2 , \theta_3 )\notag
\end{equation}
\begin{equation}
P_4 = P(2\pi-\theta_2 , 2\pi-\theta_3)\notag
\end{equation}

We also have points $P_{1,0}=P(\theta_2 (0),\theta_3 (0))$ and $P_{2,0}=P(\theta_2 (0),2\pi-\theta_3 (0))$, and $P_{i,0}$ are defined similarly for $i=3,4$.  The points $P_{i,\pi}$ are similarly defined.  

Since $\theta_3 (\pi) = \pi$, the points $P_{1,\pi}$ and $P_{2,\pi}$ are equal.  Likewise, $P_{3,\pi}$ and $P_{4,\pi}$ are equal.  But $\mct$ is not collapsible, so $\theta_2 (\pi) \neq \pi$.  So $P_{1,\pi}\neq P_{3,\pi}$ and therefore $p^{-1}(\pi)$ consists of two points.  

Since neither $\mct$ nor $\mch$ is foldable, neither $\theta_2 (0)$ nor $\theta_3 (0)$ is equal to 0, and so the points $P_{i,0}$ are distinct.  Therefore $p^{-1}(0)$ consists of four points.

We can make a continuous path of deformations between the points $P_1$ and $P_2$ via $P_{1,\pi}=p_{2,\pi}$.  Therefore, the points $P_1$ and $P_2$ are in the same path component.  Similarly, the points $P_3$ and $P_4$ are in the same path component. 

But $\theta_2$ cannot equal 0 or $\pi$, so there is no path of configurations from $P_1$ to $P_3$ and from $P_2$ to $P_4$.  So there are two path components in the moduli space.

\end{proof}

\begin{figure}[ht]
\scalebox{.60}{ \includegraphics{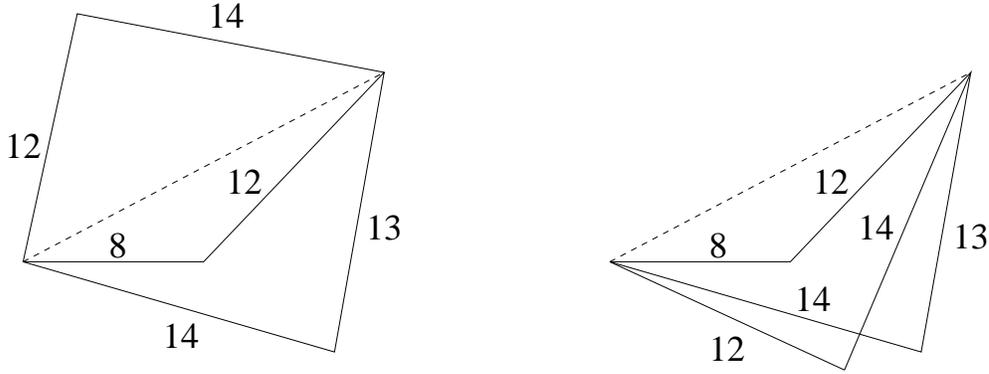}}
\caption{We cannot change the $\mct$ configuration in the first figure (above the dashed line) into the $\mct$ configuration in the second figure (below the dashed line) without dismantling the linkage.}
\label{Fi:diffcomp}
\end{figure}

\begin{figure}[ht]
\scalebox{0.60}{ \includegraphics{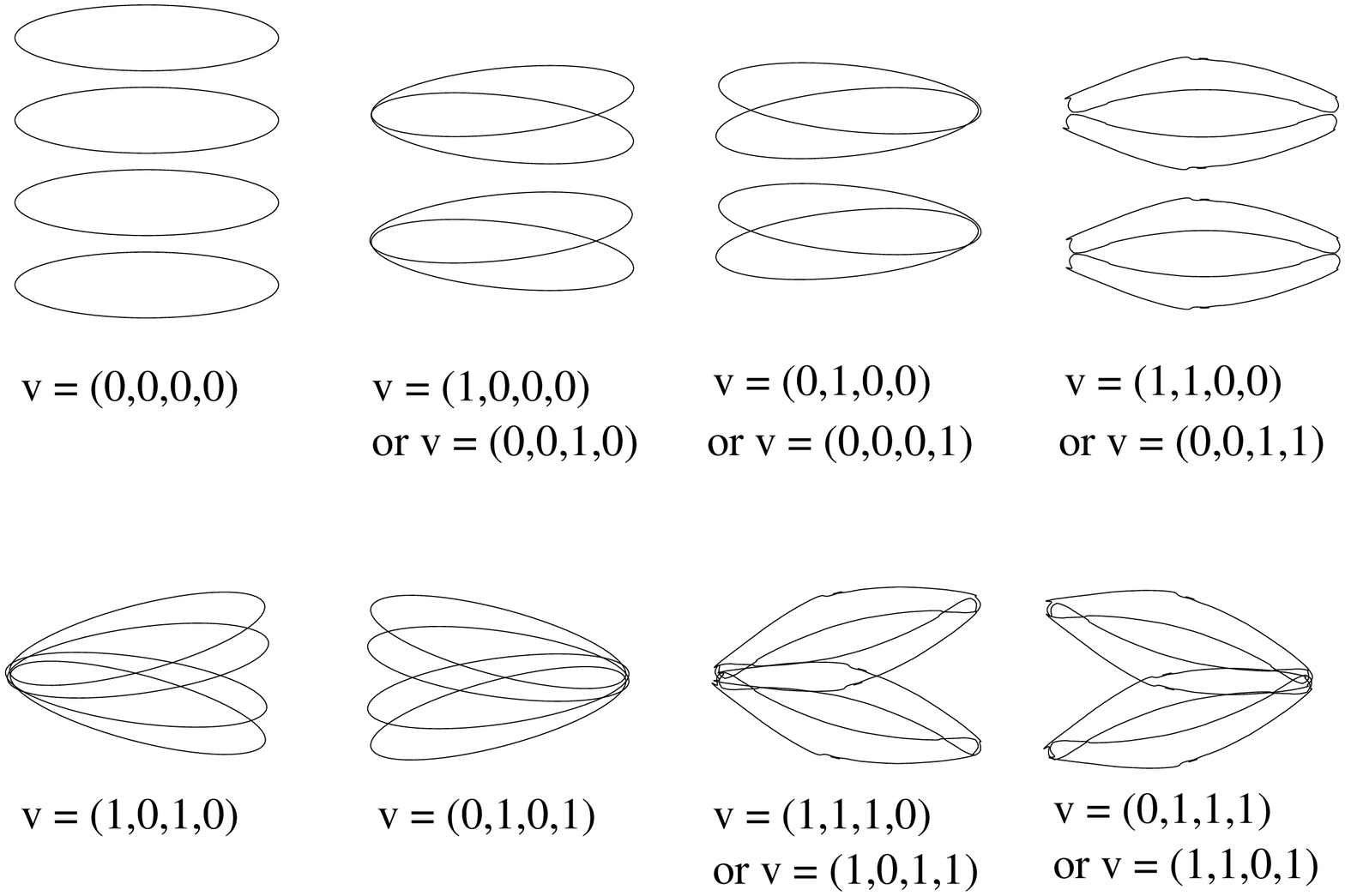}}
\caption{Case ({\bf{A}}) Moduli Spaces of Multiquadrilateral Linkages, Part 1}
\label{Fi:vcases}
\end{figure}

\begin{figure}[ht]
\scalebox{0.60}{ \includegraphics{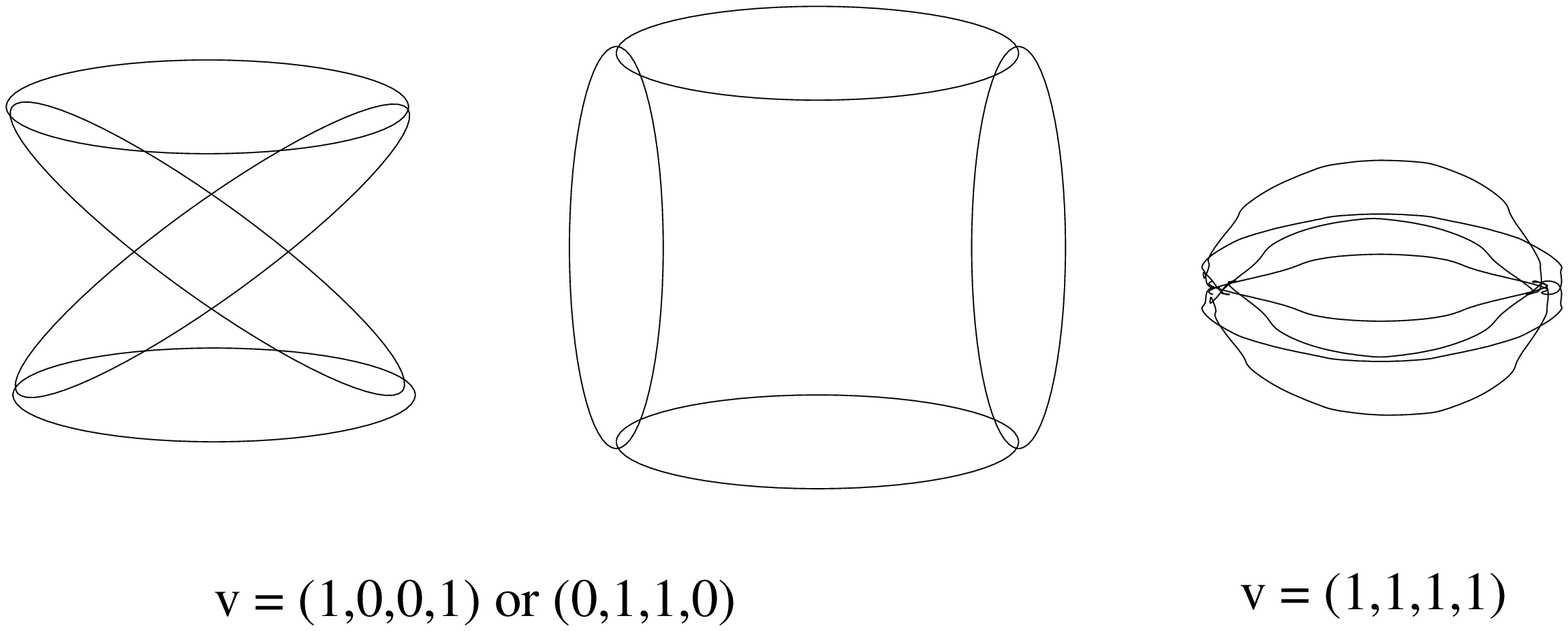}}
\caption{Case ({\bf{A}}) Moduli Spaces of Multiquadrilateral Linkages, Part 2}
\label{Fi:chaincase}
\end{figure}

We next consider case $(\mathbf{B})$.  Suppose $im(p) = [\theta_{min},\theta_{max}] \subseteq{(0,2\pi)}$.  We define a new vector $\vv=\vev$ of zeros and ones as follows.  We let 
\[
v_1 =
\begin{cases}
1 &\text{if $\mct$ is collapsible;}\\
0 &\text{if $\mct$ is not collapsible,}
\end{cases} 
\]

\[
v_2 =
\begin{cases}
1 &\text{if $\mct$ is foldable at $\theta_{min}$ and $\theta_{max}$;}\\
0 &\text{if $\mct$ is not foldable at $\theta_{min}$ or $\theta_{max}$,}
\end{cases} 
\]

\[
v_3 =
\begin{cases}
1 &\text{if $\mch$ is collapsible;}\\
0 &\text{if $\mch$ is not collapsible,}
\end{cases} 
\]
and
\[
v_4 =
\begin{cases}
1 &\text{if $\mch$ is foldable at $\theta_{min}$ and $\theta_{max}$;}\\
0 &\text{if $\mch$ is not foldable at $\theta_{min}$ or $\theta_{max}$.}
\end{cases} 
\]
We also define the pairs of $\vev$ as $({v_1},{v_2})$ and $({v_3},{v_4})$.

\begin{theorem} 
\label{T:minmaxangle}(case $(\mathbf{B})$)  Suppose that $(\mco,\mct,\mch)$ is generic and that $\pi$, but not $0$, is a possible angle for $\theta_1$.  Then the moduli space of the multiquadrilateral linkage $\mpfoth$ is formed by two circles if $\bf{v}$ has no (1,1) pairs or two figure-8's if $\bf{v}$ has at least one (1,1) pair.  Furthermore, we have the following.
\begin{align}
&\text{Number of components\qquad} &&=2^{\text{\# of (0,0) pairs in $\vv$}}\notag\\
&\text{Points in $p^{-1}(\theta_{min})=p^{-1}(\theta_{max})$\qquad} &&=2^{2-\text{\# of 1's in even positions of $\vv$}}\notag\\
&\text{Points in $p^{-1}(\pi)$\qquad}&&=2^{2-\text{\# of 1's in odd positions of $\vv$}}\notag
\end{align}
\end{theorem}

\begin{proof}  We start by showing the existence of a circle or figure-8 of configurations.  Instead of mapping $\mpfoth \rightarrow S^1$, we now map $p: \mpfoth \rightarrow [\theta_{min}, \theta_{max}]$.  Suppose $\mct$ forces $\theta_1$ to lie in that interval.  Then there is a set of configurations where the angle $\theta_3$ in $\mch$ lies in $[0,\pi]$, and there is a set of configurations where the angle $\theta_3$ lies in $[\pi,2\pi]$.  Assume that $\theta_3 \in [0,\pi]$.  
Since $\mct$ is folded when $\theta_1$ equals $\theta_{min}$ or $\theta_{max}$, we see that $p^{-1}(\theta_{min})$ and $p^{-1}(\theta_{max})$ are both one point.  If $\mct$ is not collapsible, then $p^{-1}(\theta_1)$ is two points for every $\theta_1 \in (\theta_{min},\theta_{max})$.  This gives us a circle of configurations.  If $\mct$ is collapsible, then $p^{-1}(\theta_1)$ is two points for every $\theta_1 \in (\theta_{min},\theta_{max})\backslash \{{\pi}\}$, and $p^{-1}(\pi)$ is one point.  This gives us a figure-8 of configurations.

Similarly, we obtain respective configurations when $\theta_3 \in [\pi,2\pi]$.  
\end{proof}

\section{General Results}
\label{S:gen}
\subsection{Dimension and Smoothness}
\label{dimsmooth}

Under certain conditions, the moduli space of a multipolygonal linkage is a smooth manifold.  We utilize the implicit function theorem and compute a Jacobian based on equations whose zero set is the moduli space.  Although the moduli spaces are algebraic varieties, we use equations based on the angles $\theta_{i,j}$ between the edges $e_{i,j}$ and the positive $x$-axis.

We think of each $e_{i,j}$ as a complex number $d_{i,j} e^{\sqrt{-1} \theta_{i,j}}$.  The sum of real parts of the complex numbers obtained from edges of $\mco$ must equal the sum of real parts of complex numbers obtained from edges of $\mct$ and of $\mch$.  Likewise, we can say a similar statement about the imaginary parts.
 
One case in which the moduli space is a smooth manifold is when the moduli space is described by the following equations in the angle variables.

\begin{align}
a + \cos\thto +\ldots + \cos\thnoo =&\cos\thot + \ldots +\cos\thntt\notag\\
\sin\thto + \ldots + \sin\thnoo =&\sin\thot +\ldots + \sin\thntt \notag\\
a + \cos\thto + \ldots + \cos\thnoo =& \cos\thoh + \ldots +\cos\thnhh \notag\\
\sin \thto + \ldots + \sin\thnoo =& \sin\thoh + \ldots + \sin\thnhh \notag
\end{align}
Now we create functions whose zero set is the moduli space.

\begin{align}
{f_1}= &\ a + \cos\thto +\ldots + \cos\thnoo - \cos\thot - \ldots -\cos\thntt\notag\\
{f_2}= &\ \sin\thto + \ldots + \sin\thnoo - \sin\thot - \ldots - \sin\thntt \notag\\
{f_3}= &\ a + \cos\thto + \ldots + \cos\thnoo - \cos\thoh - \ldots - \cos\thnhh \notag\\
{f_4}= &\ \sin \thto + \ldots + \sin\thnoo - \sin\thoh - \ldots - \sin\thnhh \notag
\end{align}

\begin{theorem}
\label{T:gensmooth}
The moduli space $M(a,1^{<{n_1}-1>};1^{<n_2>};1^{<n_3>})$ is a smooth manifold of dimension $n_1 + n_2 + n_3 - 5$ when the following conditions are met.
\begin{enumerate}
\item The numbers $n_2$ and $n_3$ do not have the same parity\label{I:noparity}
\item The length $a$ lies in the interval $(0,b)\setminus \mathbb{Z}$, where $b=min\{{n_2},{n_3}\}$\label{I:existnofold}
\end{enumerate}

\end{theorem}

\begin{proof}

The condition ~(\ref{I:existnofold}) above assures us that the linkage is actually possible since $a<b$ and that the first free linkage has non-integer length between $x_{1,1}$ and $x_{{n_1}+1,1}$ if it lies completely in a straight line.  We assume $n_2$ is odd and $n_3$ is even.  The proof when $n_2$ is even and $n_3$ is odd goes the same way.

Here is a heuristic dimension count.  We have $n_1 - 1$ degrees of freedom for the linkage with the $\thio$ angles, as the first edge lies on the positive $x$-axis.  But we only have $n_2 - 2$ degrees of freedom for the linkage with the $\thit$ angles.  This is because, once we have chosen the first $n_2 - 2$ angles, we must have the last two edges connect up to the terminal point of the first linkage.  We cannot do this if $|x_{n_2 -2,2}-x_{n_2 , 2}| > 2$, and we can do this in only one or two ways if $|x_{n_2 -2,2}-x_{n_2 , 2}| \leq 2$.  A similar argument shows that we only have $n_3 - 2$ degrees of freedom for the linkage with the $\thih$ angles.  The total number of free angles is $n_1 + n_2 + n_3 -5$, which equals the difference in the number of equations and unknowns above.

The matrix of the Jacobian can be written in block form 

\[
\begin{pmatrix}
A & B_1 & 0\\
0 & B_2 & C
\end{pmatrix}.
\]
The entries of $A$ come from taking partials of the first two equations with respect to the $\thit$ angles.  The entries of $B_1$ come from taking partials of the first two equations with respect to the $\thio$ variables.  The entries of $B_2$ come from taking partials of the last two equations with respect to the $\thio$ variables.  The entries of $C$ come from taking partials of the last two equations with respect to the $\thih$ variables.

\[
A=\begin{pmatrix}
\sin\thot & \cdots & \sin\thntt\\
-\cos\thot & \cdots & -\cos\thntt
\end{pmatrix}
\]

\[
{B_1}=\begin{pmatrix}
-\sin\thto & \cdots & -\sin\thnoo\\
\cos\thto & \cdots & \cos\thnoo
\end{pmatrix}
\]

\[
{B_2} = \begin{pmatrix}
-\sin\thto & \cdots & -\sin\thnoo\\
\cos\thto & \cdots & \cos\thnoo
\end{pmatrix}
\]

\[
C=\begin{pmatrix}
\sin\thoh & \cdots & \sin\thnhh\\
-\cos\thoh & \cdots & -\cos\thnhh
\end{pmatrix}
\]

We now state some lemmas whose proofs are easy.

\begin{lemma}
\label{L:afullrank}
The following are equivalent.
\begin{enumerate}
\item The rank of $A$ is 2.
\item The edges of $\mct$ do not lie in a single straight line.
\item The Jacobian is column equivalent to 
\[
\begin{pmatrix}
A & 0 & 0\\
0 & B_2 & C
\end{pmatrix}.
\]
\end{enumerate}
\end{lemma}

\begin{lemma}
\label{L:cfullrank}
The following are equivalent.
\begin{enumerate}
\item The rank of $C$ is 2.
\item The edges of $\mch$ do not lie in a single straight line.
\item The Jacobian is column equivalent to 
\[
\begin{pmatrix}
A & B_1 & 0\\
0 & 0 & C
\end{pmatrix}.
\]
\end{enumerate}

\end{lemma}

\begin{lemma}
\label{L:aline}
The following are equivalent.
\begin{enumerate}
\item The rank of $A$ is one and the rank of $[A|{B_1}]$ is 2.
\item All the edges of $\mct$ lie in a straight line which is different from the $x$-axis.
\end{enumerate}
\end{lemma}

\begin{lemma}
\label{L:cline}
The following are equivalent.
\begin{enumerate}
\item The rank of $C$ is one and the rank of $[{B_2}|C]$ is 2.
\item All the edges of $\mch$ lie in a straight line which is different from the $x$-axis.
\end{enumerate}
\end{lemma}

Using these lemmas, we can show that if $a$ is not an integer, if $n_2$ is odd, and if $n_3$ is even, then the Jacobian has full rank.  

If the edges of $\mct$ in a configuration of the linkage do not lie in a straight line, and if the edges of $\mch$ in that configuration do not lie in a straight line, then both $A$ and $C$ have rank 2.  Therefore the Jacobian has full rank.

If all of the edges of $\mct$ lie on the $x$-axis, then $A$ has rank 1 and the distance between the initial and terminal points is odd.  Since $\mch$ has an even number of edges, it cannot possibly lie in a straight line.  This forces $C$ to have rank 2.  Therefore the Jacobian is column equivalent to 
\[
\begin{pmatrix}
A & B_1 & 0\\
0 & 0 & C
\end{pmatrix}.
\]

Because $\mco$ has an edge $a$ which is not integer, the remaining edges of $\mco$ cannot all lie in the same straight line as $\mct$.  Therefore $B_1$ has a column which is not a scalar multiple of any column in $A$.  Therefore $[A|{B_1}]$ has rank 2, and the Jacobian has full rank.  A similar argument holds if we assume all of the edges of $\mch$ lie on the $x$-axis.

If all the edges of $\mct$ lie in a straight line which is not the $x$-axis, then the rank of $[A|{B_1}]$ is 2.  Since the distance between the initial and terminal points is odd, the edges of $\mch$ cannot all lie in that straight line.  Therefore $C$ has rank 2, and the Jacobian is column equivalent to
\[
\begin{pmatrix}
A & B_1 & 0\\
0 & 0 & C
\end{pmatrix}.
\]
Therefore the Jacobian has full rank.  A similar argument holds if we assume all of the edges of $\mch$ lies in a straight line which is not the $x$-axis.

\end{proof}

Now we determine a nice class of smooth manifold moduli spaces described by the following equations.

\begin{align}
{f_1}= &a + b \cos\thto - {d_{1,2}}\cos\thot - \ldots -{d_{{n_2},2}}\cos\thntt \notag\\
{f_2}= &b \sin\thto - {d_{1,2}}\sin\thot - \ldots - {d_{{n_2},2}}\sin\thntt \notag\\
{f_3}= &a + b \cos\thto - {d_{1,3}}\cos\thoh - \ldots - {d_{{n_3},3}}\cos\thnhh \notag\\
{f_4}= &b \sin \thto - {d_{1,3}}\sin\thoh - \ldots - {d_{{n_3},3}}\sin\thnhh \notag
\end{align}

\begin{theorem}
\label{T:nointsmooth}
The moduli space $M(a,b;{\bf{d_2}};{\bf{d_3}})$ is a smooth manifold of dimension $n_2 + n_3 - 3$ when $a+b<min\{ \sum_{j=i}^{n_2}{d_{j,2}},\sum_{j=1}^{n_3}{d_{j,3}} \}$ and one of the following conditions are met.
\begin{enumerate}
\item The sets ${D_2} = \{ \pm {d_{1,2}}\pm\ldots\pm {d_{{n_2},2}} \}$, ${D_3}=\{ \pm {d_{1,3}}\pm\ldots\pm {d_{{n_3},3}} \}$ and \\$D = \{ |a-b|,a+b \}$ are pairwise disjoint.\label{I:noparityanb}
\item The sets $D_2$ and $D_3$ have nonempty intersection but the sets \\$[|a-b|,a+b]\cap ({D_2}\cap {D_3})$ and $D \cap ({D_2} \cup {D_3})$ are empty.
\end{enumerate}
\end{theorem}
\begin{proof}

We build a matrix similar to that in the proof of Theorem \ref{T:gensmooth}.  Again we start by writing functions whose zero set is the moduli space.

\begin{align}
{f_1}= &a + b \cos\thto - {d_{1,2}}\cos\thot - \ldots -{d_{{n_2},2}}\cos\thntt \notag\\
{f_2}= &b \sin\thto - {d_{1,2}}\sin\thot - \ldots - {d_{{n_2},2}}\sin\thntt \notag\\
{f_3}= &a + b \cos\thto - {d_{1,3}}\cos\thoh - \ldots - {d_{{n_3},3}}\cos\thnhh \notag\\
{f_4}= &b \sin \thto - {d_{1,3}}\sin\thoh - \ldots - {d_{{n_3},3}}\sin\thnhh \notag
\end{align}

The dimension count goes similarly as in the proof of Theorem \ref{T:gensmooth}.

The Jacobian has the block form 
\[
\begin{pmatrix}
A & B_1 & 0\\
0 & B_2 & C
\end{pmatrix}
\]  
where $A$ and $C$ are the same as before.  Now $B_1$ is just the single column
$$(-b \sin\thto, b \cos\thto )^T$$
and $B_2$ is the single column
$$(-b \sin\thto,b \cos\thto)^T.$$
So we want to show that, under the conditions of \ref{T:nointsmooth}, the Jacobian has full rank.  We will use Lemmas \ref{L:afullrank}, \ref{L:cfullrank}, \ref{L:aline}, and \ref{L:cline}

Suppose that the sets $D$, $D_2$, and $D_3$ are pairwise disjoint.  In particular, $|a-b|$ and $a+b$ do not lie in $D_2$ and $D_3$.  Then the edges of $\mct$ and of $\mch$ cannot possibly all lie on the $x$-axis.  So if the edges of either $\mct$ or of $\mch$ all lie in a single straight line, then that line is different from the $x$-axis.  If neither $\mct$ nor $\mch$ have all edges in a straight line, then the Jacobian has full rank.

So suppose the edges of $\mct$ lie in a straight line which is different from the $x$-axis.  The distance between the initial and terminal points is equal to some $\pm {d_{1,2}} \pm \ldots \pm {d_{{n_2},2}}$, which does not equal any point in $D_3$.  So the edges of $\mch$ do not all lie in the same straight line.  So $C$ has rank 2, $[A|{B_1}]$ has rank 2, and the Jacobian is column equivalent to 
\[
\begin{pmatrix}
A & B_1 & 0\\
0 & 0 & C
\end{pmatrix}.
\]
Thus the Jacobian has full rank.  A similar argument follows if the edges of $\mch$ lie in a straight line which is different from the $x$-axis.

Now suppose that the sets $D_2$ and $D_3$ have nonempty intersection, but that $[|a-b|,a+b]$ has empty intersection with ${D_2}\cap{D_3}$ and that $D$ has empty intersection with ${D_2}\cup {D_3}$.  Now regardless of the angle between $a$ and $b$, the distance between the initial and terminal points will not equal any number that can be written as $\pm {d_{1,2}} \pm \ldots \pm {d_{{n_2},2}}$ and as $\pm {d_{1,3}} \pm \ldots \pm {d_{{n_3},3}}$. 

If it can be written as $\pm {d_{1,2}} \pm \ldots \pm {d_{{n_2},2}}$ but not as $\pm {d_{1,3}} \pm \ldots \pm {d_{{n_3},3}}$, then $\mct$ lies in a straight line.  But the second empty intersection condition means that this line must be different from the $x$-axis.  Then $[A|{B_1}]$ has rank 2 and $C$ has rank 2, so the Jacobian has full rank.  A similar argument follows if the distance between the initial and terminal points can be written as $\pm {d_{1,3}} \pm \ldots \pm {d_{{n_3},3}}$ but not as $\pm {d_{1,2}} \pm \ldots \pm {d_{{n_2},2}}$.

\end{proof}

Here is a special case of the theorem, where all edges $d_{i,j}$ for $j=2$ and $j=3$ are equal to 1.  If $n_2$ and $n_3$ are both even, then $D_2$ and $D_3$ and their unions and intersections are all the even integers.  If $n_2$ and $n_3$ are both even, then $D_2$ and $D_3$ and their unions and intersections are all the odd integers.  But if $n_2$ and $n_3$ do not have the same parity, then $D_2 \cap D_3 $ is empty, and $D_2 \cup D_3 $ is $\Z$.

\begin{corollary}
\label{C:nointsmooth}
The moduli space $M(a,b;1^{<n_2>};1^{<n_3>})$ is a smooth manifold of dimension $n_2 + n_3 - 3$ when $a+b<min\{ {n_2},{n_3} \}$ and one of the following conditions are met.
\begin{enumerate}
\item The numbers $n_2$ and $n_3$ do not have the same parity, $a+b\notin{\Z}$, and $|a-b|\notin{\Z}$.\label{I:noparityanb}
\item The numbers $n_2$ and $n_3$ are both odd, and $[|a-b|,a+b]$ has empty intersection with the odd integers.
\item The numbers $n_2$ and $n_3$ are both even, and $[|a-b|,a+b]$ has empty intersection with the even integers (note that this condition disallows $a=b$).
\end{enumerate}
\end{corollary}
\begin{proof}

\end{proof}

\subsection{Moduli Space as a Fibered Product}
\label{oneshared}

\begin{definition}
\label{D:fiberproduct}
Suppose $Y_1$, $Y_2$, and $Z$ are spaces and suppose $q_1 : {Y_1} \rightarrow Z$, and $q_2 : {Y_2} \rightarrow Z$ are maps.  The subspace of ${Y_1}\times{Y_2}$ defined by 
\begin{equation}
\{(y_1 , y_2 ) \in {Y_1}\times {Y_2}\ |\ {q_1}(y_1) = {q_2}(y_2)\}\notag
\end{equation}
is called the fibered product of ${Y_1}$ and ${Y_2}$ relative to $Z$, and is denoted ${Y_1}\times_Z {Y_2}$.
\end{definition}

Now suppose $Y_1 = M(\mco,\mct)$ and $Y_2 = M(\mco,\mch)$ and $Z=M(\mco)$.  The first two spaces are moduli spaces of polygonal linkages and the last space is the moduli space of one free linkage.  We have the following homeomorphism result.

\begin{theorem} 
\label{T:fiberproduct}
The space $\mpfoth$ is homeomorphic to \newline $M(\mco,\mct)\times_{M(\mco)} M(\mco,\mch)$.  
\end{theorem}

Note:  For this theorem we could possibly have an empty moduli space.  So $\emptyset \times_Z X = \emptyset$ for every space $X$.

\begin{proof}
We will construct a continuous function $$f: \mpfoth  \rightarrow \fpr$$ and a continuous inverse function $f^{-1}$.  

\begin{figure}[ht]
\scalebox{0.60}{ \includegraphics{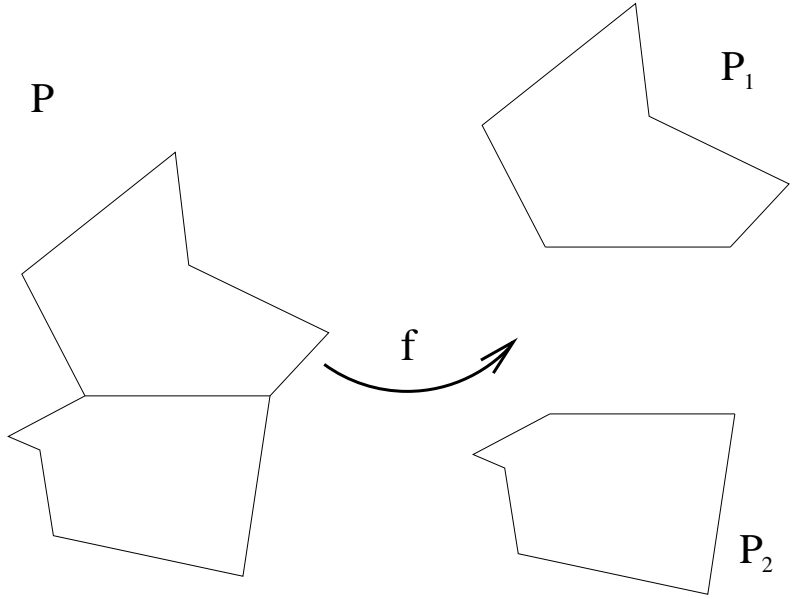}}
\caption{The effect of the $f$ map when $\mco$ has one edge}
\label{Fi:splitlink}
\end{figure}

Suppose $P$ is an embedding of $\mlink$.  We use the symbol $\prest$ to denote the embedding of $(\mco,\mct)$ after removing all edges of the $\mc$ linkage from that embedding. 

Also suppose that $P_1$ and $P_2$ are embeddings of $(\mco,\mct)$ and $(\mco,\mch)$ respectively such that $\mco$ is embedded the same way for each.   We use the notation $P_1 \sharp P_2$ to represent the embedding of $\mlink$ where the common linkage is identified.

We define $f(P)=(\prestoh,\presttt)$ and $f^{-1}(P_1 , P_2) = P_1 \sharp P_2$.  The composition $f(f^{-1})$ takes a pair of embeddings of polygonal linkages with identical configurations of $\mco$, identifies the edges of $\mco$, and then maps this embedding to a 2-tuple of embeddings, one with the $\mct$ removed and one with the $\mch$ removed.  It is clear to see that this operation is the identity on $\fpr$.  

The composition $f^{-1}(f)$ takes an embedding of the multipolygonal linkage, then maps this embedding to a 2-tuple of embeddings, one with the $\mct$ removed and one with the $\mch$ removed.  This pair of embeddings has the same configuration of the $\mco$ free linkage, so it lies in the fibered product.  So the pair is then mapped to one multipolygonal linkage embedding where the common edges are identified.  This is clearly the identity on $\mpfoth$.

\end{proof}

Notice that if $\mco$ is formed by one edge, then $M(\mco)$ is just a point, and the fibered product $\fpr$ becomes the usual cross product $M(\mco,\mct)\times M(\mco,\mch)$.

\subsection{Disjoint Union of Polygonal Moduli Spaces}
\label{disjunion}

\begin{definition}
\label{D:lengthrange}
The length range of a free linkage with vertices $\{x_i \}_{i=1}^n$ is the closed interval $[{l_1},{l_2}]$, where $l_1$ is the minimum value of $|x_{n+1}-{x_1}|$ over all points in the moduli space of $\mc$, and ${l_2}=\sum_{i=1}^n |x_{i+1}-{x_i}|$.
\end{definition}

\begin{figure}[ht]
\scalebox{0.60}{ \includegraphics{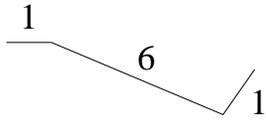}}
\caption{The length $d$ lies in the intersection of the length range intervals $[0,22]$ and $[4,8]$}
\label{Fi:lengthrange}
\end{figure}

\begin{proposition} 
\label{Pr:disjoint}
Suppose the free linkage $\mathcal{F}_1$ has length range $[{l_{11}},{l_{12}}]$ and suppose the free linkage $\mathcal{F}_2$ has length range $[{l_{21}},{l_{22}}]$.  Furthermore, suppose $\mathcal{F}_3$ consists of two edges of lengths $a$ and $b$, where $|a-b|<max\{ {l_{11}},{l_{21}} \}$ and $a+b>min\{ {l_{12}},{l_{22}} \}$.  Then $M(\mathcal{F}_1 , \mathcal{F}_2 , \mathcal{F}_3 )$ equals the disjoint union of two copies of $M( \mathcal{F}_1 , \mathcal{F}_2 )$.
\end{proposition}

\begin{figure}[ht]
\scalebox{0.60}{ \includegraphics{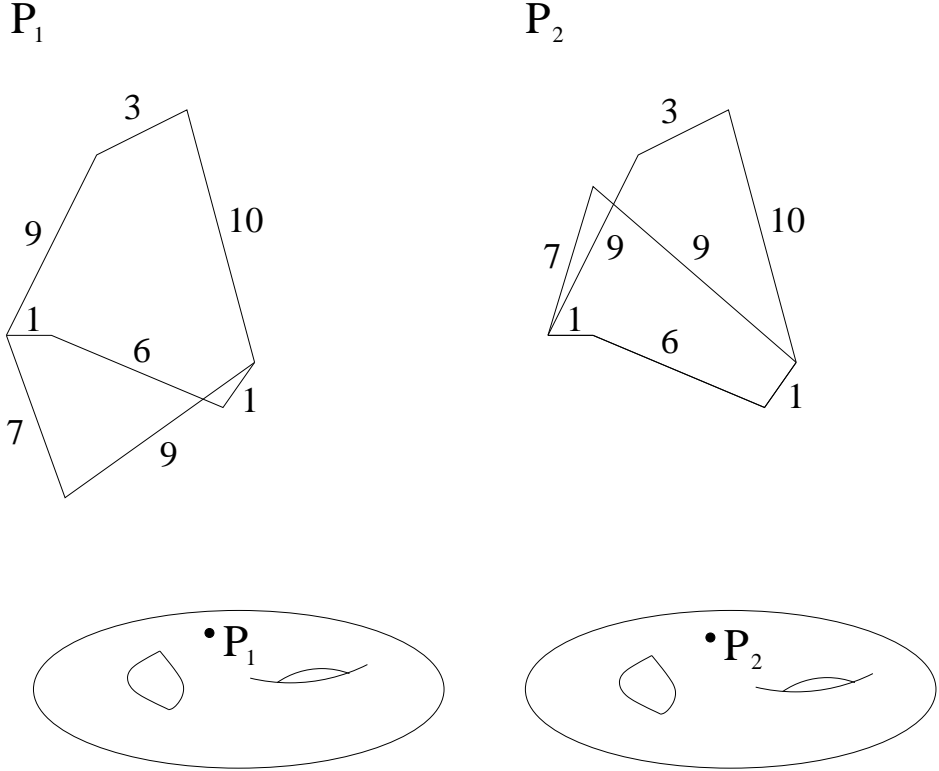}}
\caption{The added linkage with lengths 7 and 9 has length range $[2,16]$, whose interior contains the interval of the intersection of the length ranges of the other two free linkages.  There is no path of configurations from point $P_1$ to $P_2$.}
\label{Fi:pointsintwo}
\end{figure}

\begin{proof}

Let $P$ be any point in $M( \mathcal{F}_1 , \mathcal{F}_2 )$.  We will construct two points in \newline$M(\mathcal{F}_1 , \mathcal{F}_2 , \mathcal{F}_3 )$ by attaching $\mathcal{F}_3$.  The conditions $|a-b|<max\{ {l_{11}},{l_{21}} \}$ and $a+b>min\{ {l_{12}},{l_{22}} \}$ ensure that we can attach the linkage so that the initial point agrees with the identified initial points of $\mathcal{F}_1$ and $\mathcal{F}_2$ and the terminal point agrees with the identified terminal points of $\mathcal{F}_1$ and $\mathcal{F}_2$.  There are two ways to make this attachment, one in which the angle between $a$ and $b$ is in the interval $(0,\pi)$, and one in which the angle is in the interval $(\pi , 2 \pi)$.  Now because of the conditions above, there are no configurations where the angle between $a$ and $b$ is $0$ or $\pi$.  Since the angle between $a$ and $b$ has to change continuously as we deform any linkage in the moduli space, there is no path that we can take to get from one point to the other. Refer to figures \ref{Fi:lengthrange} and \ref{Fi:pointsintwo}.

\end{proof}

\subsection{Disjoint Union Modulo a Subspace}
\label{disjointmod}

Now if $|a-b|$ can equal $max\{ {l_{11}},{l_{21}} \}$ or if $a+b$ can equal $min\{ {l_{12}},{l_{22}} \}$ (or both), then there is a path between the points described above.  This is because the angle between $a$ and $b$ can equal $0$ or $\pi$ (or both).  With this observation we have the following propositions.

\begin{proposition}
\label{P:disjointmodulo1}
Suppose we have the same conditions as in \ref{Pr:disjoint} except that \newline $|a-b|=max\{ {l_{11}},{l_{21}} \} > 0$.  Then the moduli space $M(\mathcal{F}_1 , \mathcal{F}_2 , \mathcal{F}_3 )$ is a disjoint union of two copies of $\mpfot$, modulo the subspace ${M_0}^{\prime}$ of $M(\mathcal{F}_1 , \mathcal{F}_2)$ where the angle between $a$ and $b$ is $0$.
\end{proposition}

Similarly we have the following.

\begin{proposition}
\label{P:disjointmodulo2}
Suppose we have the same conditions as in \ref{Pr:disjoint} except that $a+b=min\{ {l_{12}},{l_{22}} \}$.  Then the moduli space $M(\mathcal{F}_1 , \mathcal{F}_2 , \mathcal{F}_3 )$ is a disjoint union of two copies of $\mpfot$, modulo the subspace ${M_{\pi}}^{\prime}$ of $M(\mathcal{F}_1 , \mathcal{F}_2)$ where the angle between $a$ and $b$ is $\pi$.
\end{proposition}

Since the moduli space is the union of two algebraic subvarieties, which are in turn triangulable, we have the following Mayer-Vietoris sequence.

\begin{equation}
\ldots \rightarrow H_i ({M_{\theta}}^{\prime})  \rightarrow H_i (M(\mco,\mct)) \oplus H_i (M(\mco,\mct)) \rightarrow H_i (\mpfoth) \rightarrow \ldots \notag
\end{equation}

This gives us the following Euler characteristics.

\begin{equation}
\chi(\mpfoth)=2 \chi(\mpfot)-\chi({M_{\theta}}^{\prime})\notag
\end{equation}

\subsection{Connectedness}
\label{connect}

We recall Theorem 1 from \cite{KM1} regarding polygonal linkages.  Kapovich and Millson use $M_r$ as the notation for the moduli space of a polygonal linkage with a vector $r=({r_1},\ldots,{r_n})$ of edge lengths.  They also allow some, but not all, edges to have a possible length of 0 and normalize their polygons so that the perimeter is 1.  That theorem generalizes to the following form.

\begin{theorem}
\label{T:KMnonorm}
The space $M_r$ is not connected if and only if there are three different sides $e_i$, $e_j$, and $e_k$ (with lengths $r_i$, $r_j$ and $r_k$ respectively) in the normalized polygon so that
\begin{equation}
{r_i}+{r_j}>{\frac{1}{2}}\sum_{i^{\prime}=1}^n {r_i^{\prime}},\;\;\;\; {r_j}+{r_k}>{\frac{1}{2}}\sum_{i^{\prime}=1}^n {r_i^{\prime}},\;\;\;\; {r_i}+{r_k}>{\frac{1}{2}} \sum_{i^{\prime}=1}^n {r_i^{\prime}}\notag
\end{equation}
\end{theorem}

For the theorem below, we have the following proposition.
\begin{proposition}
\label{T:keychain}
Let $\mpfoth$ be the moduli space of a multipolygonal linkage, and let $p:\mpfoth\rightarrow{S^1}$ be an onto map which sends a configuration with angle $\theta_{1,1}$ to the angle ${\theta_{1,1}}$.  Suppose the following conditions hold.

\begin{enumerate}

\item $\mpfoth$ contains a subspace homeomorphic to $S^1$ which intersects each fiber.

\item For every $\theta \in im(p)$, $p^{-1} (\theta)$ is connected.
\end{enumerate}

Then $\mpfoth$ is connected.
\end{proposition}

We use this theorem to give a nice class of multipolygonal linkages with connected moduli spaces.  For this theorem, let 
$$(\mco,\mct,\mch)=(\mc(a,b),\mct(r_1 , \ldots, r_{n_2} ) ,\mch({r_1}^{\prime},\ldots, {r_{n_3}}^{\prime}))$$ be a multipolygonal linkage where $\mct$ has length range containing $[|a-b|,a+b]$ and $\mch$ has length range containing $[|a-b|,a+b]$.  Let $P^{\prime}$ be the perimeter of the polygonal linkage $(\mc(|a-b|),\mct)$ and let $P^{\prime\prime}$ be the perimeter of the polygonal linkage $(\mc(|a-b|),\mch)$

\begin{theorem}
\label{T:connected}
Let $(\mco,\mct,\mch)$ be a multipolygonal linkage as described above.  Let $r_i$ and $r_j$ be the two longest edges in $\mct$ and let ${r_i}^{\prime}$ and ${r_j}^{\prime}$ be the two longest edges in $\mch$.  Suppose that ${r_1}+{r_j}\leq{P^{\prime}/2}$ and that ${r_i}^{\prime}+{r_j}^{\prime}\leq{P^{\prime\prime}}/2$.  Then the moduli space $\mpfoth$ is connected.
\end{theorem}

\begin{proof}
The length range conditions above assure that the space $\mpfoth$ contains at least one point in each fiber, so that we can change the angle $\theta = \theta_{1,1}$ between $a$ and $b$ and obtain a copy of $S^1$ that intersects each fiber.  We do this by starting with a configuration $P$ at $\theta=0$.  Then we find two paths $z_1$ (by perturbing $P$ so that $\theta$ increases from 0 to $\pi$ ) and $z_2$ (by perturbing $P$ so that $\theta$ decreases from $2 \pi$ to $\pi$) from $P$ to a common configuration with angle $\theta=\pi$.  The composition of the paths $z_1$ and ${z_2}^{-1}$ gives us a subspace homeomorphic to $S^1$ which intersects every fiber.

Now we show that each fiber is connected.  Let $\theta$ be any angle between $a$ and $b$, and let $d$ be the distance from the initial point to the terminal point of $\mco$.  Let $P_{d,2}$ be the perimeter of the polygonal linkage $(\mc(d),\mct)$.  Let $P_{d,3}$ be the perimeter of the polygonal linkage $(\mc(d),\mch)$.  Since $d\geq{|a-b|}$, we have that $P_{d,2}\geq{P^{\prime}}$ and that $P_{d,3}\geq{P^{\prime\prime}}$.  Therefore 

\begin{equation}
{r_i}+{r_j}\leq{P_{d,2}/2}\notag
\end{equation}
and
\begin{equation}
{r_i}^{\prime}+{r_j}^{\prime}\leq{P_{d,3}/2}.\notag
\end{equation}

The new polygonal linkages therefore do not have three long edges, regardless of how long $d$ is.  Since $M(d)$ is just a point, by Section \ref{oneshared} the fiber is now a cross product of moduli spaces of these polygonal linkages.  Each of these moduli spaces is connected.  Therefore the fiber is connected.  Since $\theta$ was arbitrary, every fiber of the map is connected, and so is the moduli space $\mpfoth$.
\end{proof}

Note that in the proof of theorem \ref{T:connected}, we needed paths $z_1$ and $z_2$ between configurations with different lengths between the initial and terminal vertices.  Each of those different lengths lie in the length range of the free linkage.  

\begin{proposition}
There exists a path of configurations from a free linkage with length $a$ between its initial and terminal vertices to a free linkage with length $b$ between its initial and terminal vertices, when the length range of the free linkage is $[c,d]$ and $c\leq a < b \leq d$.
\end{proposition}

\begin{proof}

We prove this proposition by deriving an algorithm for finding a path between the two configurations.  The idea for the algorithm is to make as many angles between $x_i$ and the $x$-axis as possible to be zero, with the rest of the free linkage ``sliding down''.  Then when we cannot make some angle between $x_i$ and the $x$-axis equal zero, we minimize it.

\vspace{0.10in}

{\bf{Step 1:}}  Initialize $\mc_1 ^{\prime} = \mc (d_2 , \ldots, d_n )$.  Start with any configuration of $\mc_1 ^{\prime}$.

\vspace{0.10in}

{\bf{Step 2:}}  Given $\mc_i ^{\prime}$, let $l_i ^{\prime} = |x_{i+1}-x_{n+1}|$ and let $\theta_i$ be the angle between the segment connecting $x_i$ and $x_{n+1}$ and the edge of length $d_i$.

\vspace{0.10in}

{\bf{Step 3:}}  If the following inequality holds, then go to Step 4.  Otherwise go to Step 5.

$$\sum_{j=1}^i d_j + l_i ^{\prime} \leq b$$

\vspace{0.10in}

{\bf{Step 4:}}  Perturb $\mc_i ^{\prime}$ so that $\theta_i$ becomes 0, but keeping all other angles rigid and sliding the rest of the linkage so that $x_{i+1}$ and $x_{n+1}$ are both on the positive $x$-axis and still $l_i ^{\prime}$ units apart.  If the inequality is strict, then increment $i$ and go to Step 2.  Otherwise we are done.

\vspace{0.10in}

{\bf{Step 5:}}  Perturb $\mc_i ^{\prime}$ so that $\theta_i$ decreases to a minimum angle where there exists a triangle with lengths $d_i$, $l_i ^ {\prime}$, and $b-\sum_{j=1}^{i-1} {d_j}$.  Keep all remaining angles rigid.

\vspace{0.10in}

This algorithm terminates because the length range is at most $d=\sum_{i=1}^n {d_i}$, and $b \leq d$.

\end{proof}

\section{Specific Results}
\label{S:spec}
\subsection{Bundle Structure}
\label{bundles}

Now we investigate moduli spaces of the form \linebreak $M(\mco,\mct,\mch)$ where $\mco=\mathcal{F}(a,b)$ is now the 2-bar linkage, $\mct=\mathcal{F}(1^{<n_1>})$, and $\mch = \mathcal{F}(1^{<n_2>})$.  Let $\theta$ be the angle from the $x$-axis to $b$, and let $\theta_1 = \pi-\theta$ be the angle from $b$ to $a$.  Consider the map $p:M(\mco,\mct,\mch)\rightarrow{S^1}$ which maps a configuration to its angle $\theta_1$.  When we form a triangle with lengths $a$ and $b$, with an angle $\theta_1$ between them, we have a unique number $d(\theta_1)$ that is the length of the third side. 

So we fix $\theta_1$ and we consider $p^{-1}(\theta_1)\subset\mpfoth$.  We have a homeomorphism $\phi: p^{-1}(\theta_1)\rightarrow M(\mathcal{F}(d(\theta_1)),\mct,\mch)$ which just rotates the length $d(\theta_1)$ segment to the positive $x$-axis.  The following observation will be used later.

\begin{proposition}
\label{P:subvariety}
The space $p^{-1}(\theta_1 )$ is an algebraic subvariety of $\mpfoth$.
\end{proposition}

We combine this with the homeomorphism between $ M(\mathcal{F}(d(\theta_1)),\mct,\mch)$ and $M(\mathcal{F}(d(\theta_1)),\mct)\times M(\mathcal{F}(d(\theta_1)),\mch)$ established in section \ref{oneshared}.

But here $M(\mathcal{F}(d(\theta_1)),\mct)$ and $M(\mathcal{F}(d(\theta_1)),\mch)$ are moduli spaces of quasi-equilateral linkages, so the know something about the moduli spaces from Theorems in \cite{Kqe}.  

\begin{theorem}  
\label{T:fiber}
Suppose the multipolygonal linkage $\mlink$ satisfies any one of the following three conditions.

\begin{enumerate}

\item The interval $[|a-b|,a+b]$ does not intersect the odd integers, and $\mct$ and $\mch$ both have an odd number of edges.

\item The interval $[|a-b|,a+b]$ does not intersect the even integrs, and $\mct$ and $\mch$ both have an even number of edges.

\item The interval $[|a-b|,a+b]$ does not intersect the integers, and $\mct$ and $\mch$ have opposite parity of number of edges.

\end{enumerate}
then the map $p:\mpfoth\rightarrow{S^1}$ is a locally trivial fibration, and therefore $\chi(\mpfoth)=0$.
\end{theorem}

We need the following theorems.

\begin{theorem}
\label {T:ranking}
Let $f_1 , \ldots , f_r $ be $C^{\infty}$ functions on $\R^n$ with $r \leq{n}$ and coordinates ${x_1},\ldots,{x_n}$.  Let $a$ be a point with ${f_1}(a) =\ldots={f_r}(a)=0$.  Let $K\subseteq\{1,\ldots,n\}$ be a subset of size $r$.  Suppose that the $r\times r$ matrix 
\begin{equation}
\left( \frac{\partial {f_i}}{\partial {x_k}} \right),\notag
\end{equation}
with $i=1,\ldots,r$ and $k\in{K}$, has rank $r$.  Then
\begin{enumerate}
\item The set $M\subseteq{\R^n}$ defined by ${f_1}=\ldots={f_r}=0$ is a smooth manifold of dimension $n-r$ in a neighborhood of $a$.\\
\item For any $j\notin{K}$, the projection $M\rightarrow \R$ given by $x_j$ is a submersion in a neighborhood of $a$.
\end{enumerate}
\end{theorem}

\begin{proof}

The implicit function theorem shows in fact that $\{ x_j \}_{j \notin K}$ is a system of local coordinates for a manifold structure on $M$ in a neighborhood of $a$.
\end{proof}

\begin{theorem}
\label{T:proper}
Let $f:X\rightarrow \R$ be a smooth map of smooth manifolds.  Assume that $f$ is a proper map (that is, $f^{-1}(K)$ is compact for any compact set $K$).  Assume $f$ is a submersion everywhere on $X$.  Then $X$ is a locally trivial fibration of smooth manifolds.
\end{theorem}

\begin{proof}  

We apply the proof of Theorem 3.1 in \cite{Mil} to the function $f$, which is a Morse function without critical points. 

\end{proof}

Recall that a smooth map $f:X\rightarrow Y$ is a submersion at $x$ if and only if the map $T_x (X) \rightarrow T_{f(x)} (Y)$ of tangent spaces is onto. 

\begin{proof}[Proof of theorem ~\ref{T:fiber}]  Let $x\in X=\mpfoth$ and let $Y=S^1$.  Then $x$ represents a configuration of a multipolygonal linkage which is described as the zero set of the following equations, where $\theta_1$ is the angle from $b$ to $a$ in $x$.
 
\begin{align*}
{f_1}=&a + b \cos\theta - \cos\thot - \ldots -\cos\thntt\\
{f_2}=&b \sin\theta - \sin\thot - \ldots - \sin\thntt \\
{f_3}=&a + b \cos\theta - \cos\thoh - \ldots - \cos\thnhh \\
{f_4}=&b \sin \theta - \sin\thoh - \ldots - \sin\thnhh \\
{f_5}=&\theta + \theta_1-\pi
\end{align*}

We need to show that the Jacobian of this set of equations has full rank.  The Jacobian has the form 

\[
\begin{pmatrix}
A & B_1 & 0\\
0 & B_2 & C\\
0 & 1 & 0
\end{pmatrix}
\]
where $A$, $B_1$, $B_2$, and $C$ are defined as before. 

If the conditions of $[|a-b|,a+b]$ not intersecting integers of the proper parity above are satisfied, then neither $\mct$ nor $\mch$ can lie in a single straight line.  This will force both $A$ and $C$ to have rank 2, and the last row will contribute a 5th linearly independent row.  Therefore the Jacobian has full rank.  The first result then follows from Theorems \ref{T:ranking} and \ref{T:proper}.

To show that the Euler characteristic of $\mpfoth$ is zero, we use the fact that $\chi{(\mpfoth)}= \chi(F) \chi(S^1)$, where $F$ is the fiber over every point.  Since $\chi(S^1)=0$, we see that $\chi(\mpfoth)=0$ also.
\end{proof}

We can say more here as well.  

\begin{proposition} Consider the multi-polygonal linkage $(\mco(a,b),\mct(1^{<n_2>}),\mct(1^{<n_3>}))$, where $[|a-b|,a+b]\cup \Z$ is empty.  Then the moduli space $\mpfoth$ of the multi-polygonal linkage is homeomorphic to $M(\mco,\mct)\times M(\mco,\mch) \times S^1$.
\end{proposition}

This is a special case of the more general result.

\begin{proposition} 
\label{P:nofold}
Consider the multi-polygonal linkage 
$$(\mco(a,b),\mct(d_{2,1},\ldots,d_{2,{n_2}}),\mch(d_{3,1},\ldots,d_{3,{n_3}}))$$.  
\begin{enumerate}
\item Suppose that the interval $[|a-b|,a+b]$ has empty intersection with the sets 
$$\{\pm{d_{2,1}}\pm\cdots\pm{d_{2,{n_2}}}\}$$ and
$$\{\pm{d_{3,1}}\pm\cdots\pm{d_{3,{n_3}}}\}.$$  \label{I:nofold}
\end{enumerate}
Then the moduli space of $(\mco,\mct,\mch)$ is homeomorphic to $M(\mco,\mct)\times M(\mco,\mch) \times S^1$.
\end{proposition}

\begin{proof}
We use the following result from Kapovich and Millson.

\begin{theorem}[Corollary 15] Suppose that $r = ({r_1},\ldots,{r_n}) \in int(D_n)$ does not lie on any wall, and let $\hat{r}=({r_1},\ldots,{r_n},\epsilon)$.  Then for sufficiently small positive $\epsilon$ the space $M_{\hat{r}}$ is diffeomorphic to $M_r \times S^1$.
\end{theorem}

The vector $r$ lies on a wall if and only if there exists a configuration of the polygonal linkage where every edge lies in a single line.  This condition is equivalent to the existence of a configuration such that all but one edge lie in a single line.  In terms of edge lengths, these conditions are equivalent to the condition that $r_1$ does not equal any point in 
$$\{\pm{r_2}\pm\cdots\pm{r_n}\}.$$
The sufficiently small $\epsilon$ would then satisfy $[{r_1}-\epsilon,{r_1}+\epsilon]$ having empty intersection with
$$\{\pm{r_2}\pm\cdots\pm{r_n}\}.$$

Consider $\mco(a,b)$ and the distance $d$ between its initial and terminal points.  Let $p:\mpfoth\rightarrow S^1$ by sending a configuration with angle $\theta$ to $\theta$.  Because of condition \ref{I:nofold} in \ref{P:nofold}, none of the fibers of $p$, which are cross products of moduli spaces of polygonal linkages, ever fold up into a straight line.  Therefore they have the same topological type regardless of $\theta$ (and $d$) because they do not meet any wall.

\end{proof}

\subsection{When Fibers Change Topological Type}
\label{changetype}

In general, the fibers over different angles will not all be homeomorphic.  There can be degenerate fibers where the topology changes, and the fibers of points in open intervals between the degenerate fibers will be homeomorphic to one another.

The theorem \ref{C:nointsmooth} will now be used in a more general setting.  We recall the following variant of the Poincar\'{e} Duality Theorem from Spanier. See \cite [pages~297]{Spa}

\begin{theorem}
\label{T:poincare}
Let $(X,A)$ be a compact relative $n$-manifold such that $X-A$ is orientable over $R$.  For all $q$ and $R$-modules $G$ there is an isomorphism 

\begin{equation}
H_q (X-A;G) = {\bar{H}}^{n-q} (X,A;G),
\end{equation}
where $\bar{H}$ is the homology theory used in \cite{Spa}.
\end{theorem}

Here is a consequence of that theorem.

\begin{proposition}
\label{P:eulerchar}
Let $(X,A)$ be a relative $n$-manifold.  Suppose that $(X,A)$ is orientable over a ring $R$.  Then $\chi(X) = \chi(A) + (-1)^n \chi(X-A)$.
\end{proposition}

\begin{proof}

From the long exact homology sequence with coefficients in a ring $R$, we have 
\begin{equation}
\cdots\rightarrow H_i (A) \rightarrow H_i (X) \rightarrow H_i (X,A) \rightarrow \cdots\notag
\end{equation}

This sequence gives us the following well-known property of Euler characteristics.

\begin{equation}
\chi(X) = \chi(A) + \chi (X,A)\notag
\end{equation}

On the other hand, Theorem \ref{T:poincare} states that there is a canonical isomorphism

\begin{equation}
H_i (X-A) \simeq \bar{H}^{n-i}(X,A) \cong H^{n-i} (X,A)\notag
\end{equation}
with the second isomorphism holding because of tautness.  The tautness properties hold because $(X,A)$ admits a triangulation.  By the universal coefficient theorem, $rk(H^j (X,A)) = rk(H_j (X,A))$.

Since any manifold is $\Z/2$-orientable, we have the following for $G=\Z/2$ (or $G$ can be any $\Z/2$-module.

\begin{equation}
H_q (X-A;\Z/2) = H^{n-q} (X,A;\Z/2)\notag
=H^{n-q}(X,A;\Z/2)
\end{equation}

Now we have the following Euler characteristic computation.

\begin{align*}
\chi(X) &= \chi(A)+\chi(X,A)\\
  &=\chi(A)+\sum_{i=0}^{n} (-1)^i rk(H_i (X,A))\\
  &=\chi(A)+\sum_{i=0}^{n} (-1)^i rk(H^{i} (X,A))\\
  &=\chi(A)+\sum_{i=0}^{n} (-1)^{n-i} rk(H^{n-i} (X,A))\\
  &=\chi(A)+\sum_{i=0}^{n} (-1)^{n-i} rk(H_{i} (X-A))\\
  &=\chi(A)+(-1)^n \sum_{i=0}^{n} (-1)^{i} rk(H_{i} (X-A))\\
  &=\chi(A)+ (-1)^n \chi(X-A)
\end{align*}
\end{proof}

\begin{theorem}
\label{T:homeointerval}
Let $p:X\rightarrow S^1$ be a proper, onto map, and suppose there exists a finite set $T=\{\theta_i \}\subset{S^1}$ such that $p|_{X-p^{-1}(T)}$ is a smooth map of maximal rank.  Then the map $p$ is a locally trivial fibration on $p^{-1}(I)$ for each open interval $I$ of $S^1 - T$.
\end{theorem}

\begin{proof} The proof uses Theorem \ref{T:proper}. \end{proof}

Let $X=\mpfoth$ and let $A$ be the disjoint union of the singular fibers $\{F_i \}$.  There are finitely many such fibers.  These occur as $p^{-1}(\theta_i)$ for $\theta_i$ angles causing the distance between the initial and terminal points to be integers of the proper parity.  Then $\chi(A)=\sum\chi(F_i)$.

Then $X-A$ is the disjoint union of the inverse images of the open intervals in $S^1 - \cup \{ \theta_i \}$.  Since all fibers are homeomorphic on a fixed interval, we can contract these inverse images down to inverse images of a single point in each interval.  So $B$ is homotopic to $\sqcup{{F_j}^{\prime}}$ where the ${{F_j}^{\prime}}$ fibers are in one-to-one correspondence with the open intervals of $S^1 - \cup \{ \theta_i \}$.  Then $\chi(B) = \sum\chi({{F_j}^{\prime}})$.

As a result we obtain $\chi(\mpfoth)=\sum\chi(F_i) + (-1)^n \sum\chi({{F_i}^{\prime}})$.

\subsection{Example Where Fibers Change Topological Type}
\label{example}

\begin{figure}[ht]
\scalebox{0.60}{ \includegraphics{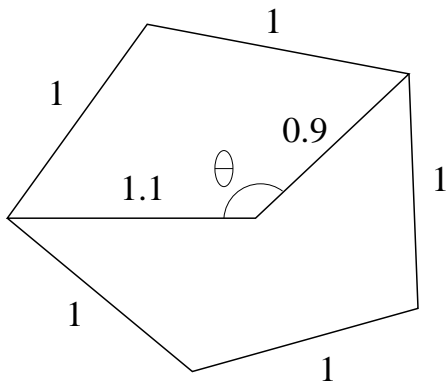}}
\caption{Sample multipolygonal linkages with three degenerate fibers}
\label{Fi:threedegenerate}
\end{figure}

For example, consider the linkage $(\mco,\mct,\mch)=(\mco(1.1,0.9),\mct(1,1,1),\mch(1,1))$ whose picture is given in figure \ref{Fi:threedegenerate}.

\begin{figure}[ht]
\scalebox{0.60}{ \includegraphics{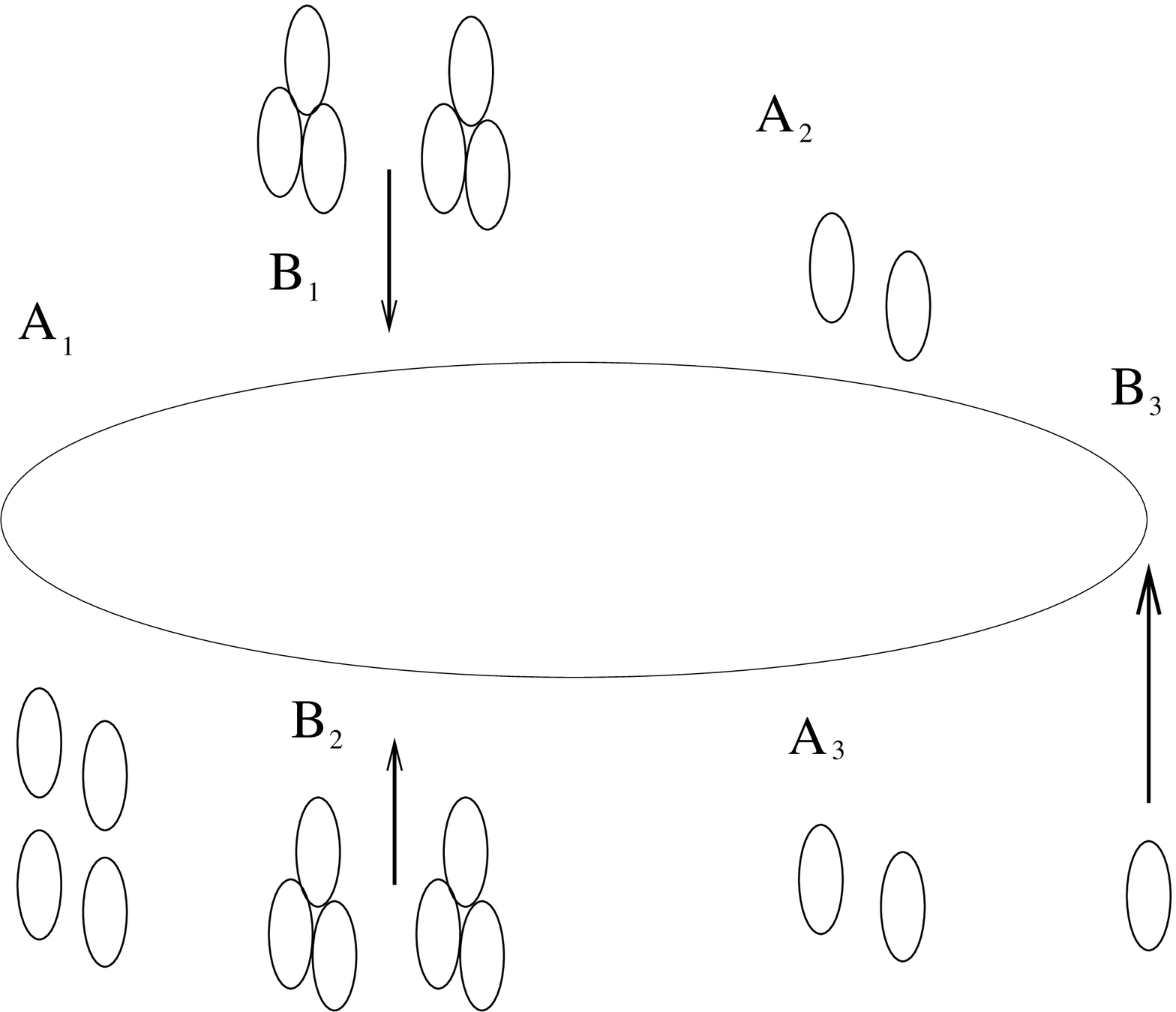}}
\caption{The fibers of the map \newline $p:M(\mco(1.1,0.9),\mct(1,1,1),\mch(1,1))\rightarrow {S^1}$}
\label{Fi:quad}
\end{figure}

When $\theta = \cos^{-1} (51/99)$ and $-\cos^{-1}(51/99)$, the distance between the initial and terminal edges equals 1.  The fibers over these two angles is the cross product of the moduli space of an equilateral quadrilateral linkage and that of a triangle.  The former moduli space is formed by three circles, any two of which meet in one point.    The latter moduli space is two points.  So the moduli space of each fiber is homeomorphic to a disjoint union of two copies of the space in Figure \ref{Fi:quad}.

\begin{figure}[ht]
\scalebox{1.00}{ \includegraphics{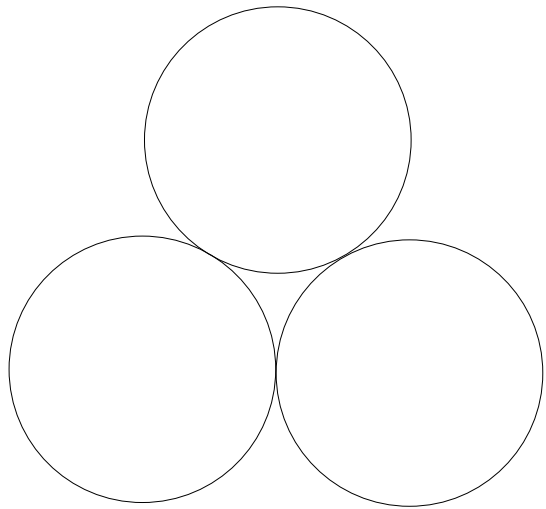}}
\caption{The moduli space of an equilateral quadrilateral linkage}
\label{Fi:equad}
\end{figure}

On the interval $(-\cos^{-1}(51/99),\cos^{-1}(51/99))$ the fiber over any angle is homeomorphic to the cross product of the moduli space of $(\mco(d),\mct)$ (where $d<1$) and that of a triangle.  The former moduli space is two circles and the latter moduli space is two points.  So the fiber over any angle in that interval is a disjoint union of four circles.

On the intervals $(-\pi,-\cos^{-1}(51/99))$ and $(\cos^{-1}(51/99),\pi)$, the fiber over any angle is homeomorphic to the cross product of the moduli space of $(\mco(d),\mct)$ (where $1<d<2$) and that of a triangle.  The former moduli space is one circle and the latter moduli space is two points.  So the fiber over any angle in that interval is a disjoint union of two circles.

Finally the fiber of $\pi$ i  homeomorphic to the cross product of the moduli space of $(\mco(2),\mct)$ and that of a triangle.  The former moduli space is one circle and the latter moduli space is one point because $\mco$ has collapsed into a straight line.  So the fiber over $\pi$ is one circle.

So $\chi(\mpfoth)$ is the sum of Euler characteristics of four circles, 2 circles, 2 circles, one circle, and four copies of the space in Figure \ref{Fi:quad}.  All except the last space have Euler characteristic 0.  The last space has Euler characteristic -3 for each copy and a total Euler characteristic of -12.  Therefore $\chi(\mpfoth)=-12$.

See the figure below.

\begin{thebibliography}{99}

\bibitem{KM2}
   Michael Kapovich and John~J. Millson,
   \emph{The Symplectic Geometry of Polygons in Euclidean Space},
   Journal of Differential Geometry \textbf{44}~(1996) 479--513.

\bibitem{KM1}
   Michael Kapovich and John~J. Millson,
   \emph{On The Moduli Space of Polygons in the Euclidean Plane},
   Journal of Differential Geometry \textbf{42}~(1995) 133--164.

\bibitem{Kmo}
   Yasuhiko Kamiyama,
   \emph{Topology of Equilateral Polygon Linkages in the Euclidean Plane Modulo Isometry Group},
   Osaka Journal of Mathematics \textbf{36}~(1999) 731--745.

\bibitem{Ksp}
   Yasuhiko Kamiyama,
   \emph{Remarks of the Topology of Spatial Polygon Spaces},
   Bulletin of the Australian Mathematical Society \textbf{58}~(1998) 373--382.

\bibitem{Kho}
   Yasuhiko Kamiyama,
   \emph{The Homology of Singular Polygon Spaces},
   Canadian Journal of Mathematics \textbf{50}~(1998) 581--594.

\bibitem{KTs}
   Yasuhiko Kamiyama and Michishige Tezuka,
   \emph{Symplectic Volume of the Moduli Space of Spatial Polygons},
   Journal of Mathematics of Kyoto University \textbf{39}~(1999) 557--575.

\bibitem{Kcn}
   Yasuhiko Kamiyama,
   \emph{Chern Numbers of the Moduli Space of Spatial Polygons},
   Kodai Mathematical Journal \textbf{23}~(2000) 380--390.

\bibitem{Kec}
   Yasuhiko Kamiyama,
   \emph{Euler Characteristics of the Moduli Space of Polygons in Higher-dimensional Euclidean Space},
   Kyushu Journal of Mathematics \textbf{54}~(2000) 333--369.

\bibitem{Kte}
   Yasuhiko Kamiyama,
   \emph{Topology of Equilateral Polygon Linkages},
   Topology Appl, \textbf{68}~(1996) 13--31.

\bibitem{Kqe}
   Yasuhiko Kamiyama and Michishige Tezuka and Tsuguyoshi Toma,
   \emph{Homology of the Configuration Spaces of Quasi-Equilateral Polygon Linkages},
   Transactions of the American Mathematical Society \textbf{350}~(1998) 4869--4896.

\bibitem{Mil}
   J. Milnor,
   \emph{Morse Theory},
   Princeton University Press , 1969.


\bibitem{Spa}
   Edwin~H. Spanier,
   \emph{Algebraic Topology},
   McGraw-Hill, Inc. , 1966.

\bibitem{Tri}
   J. Bochnak and M. Coste and M-F. Roy,
   \emph{G\'{e}om\'{e}trie Alg\'{e}brique R\'{e}elle},
   Springer Verlag , 1987.










\end{thebibliography}
\end{document}